\documentclass[a4paper,11pt,onecolumn,oneside,notitlepage,openbib,final]{article}
\pdfoutput=1
\usepackage{mathpazo}	
\usepackage{relsize}
\renewcommand{\ast}{{\mathlarger *}} 
\usepackage{float}
\usepackage{lscape}	
\usepackage[bottom]{footmisc}			
\usepackage[margin=1.2in]{geometry}
\usepackage{setspace} 
\usepackage{amsmath,amsfonts,amssymb,amsthm}
\usepackage{bbm} 
\usepackage{mathrsfs}		
\allowdisplaybreaks
\usepackage{thmtools}
\usepackage{mathtools}
\usepackage[pdflatex]{graphicx}
\usepackage{color}
\usepackage[labelformat=parens]{subfig}		
\usepackage[notcite]{showkeys}
\usepackage[pdftex,bookmarks=true,bookmarksnumbered=true,hidelinks]{hyperref}
\usepackage{url}
\usepackage{natbib}
\usepackage{enumerate}
\usepackage{alltt}	
\usepackage{cleveref} 
\usepackage[most]{tcolorbox}
\usepackage{listings} 
\usepackage{adjustbox} 
\usepackage{multirow}
\usepackage{array}
\newcolumntype{C}{>{$\displaystyle} c <{$}}

\def\tdif#1#2{d #1/d #2}

\def\tpdif#1#2{\partial #1 / \partial #2}

\def\2'{^{\prime\prime}}

\def\F{\mathbf F}

\def\R{\mathbb R}

\def\x{\mathbf x}

\def\z{\mathbf z}
\def\0{\mathbf 0}
\def\1{\mathbf 1}

\def\cX{\mathcal X}

\def\cV{\mathcal V}

\def\bpi{{\boldsymbol{\pi}}}

\def\bcV{\boldsymbol{\cV}}

\DeclareMathOperator{\bd}{bd}
 
\DeclareMathOperator{\cl}{cl}

\DeclareFontFamily{U}{mathx}{\hyphenchar\font45}
\DeclareFontShape{U}{mathx}{m}{n}{
      <5> <6> <7> <8> <9> <10>
      <10.95> <12> <14.4> <17.28> <20.74> <24.88>
      mathx10
      }{}
\DeclareSymbolFont{mathx}{U}{mathx}{m}{n}
\DeclareMathSymbol{\bigtimes}{1}{mathx}{"91}

\theoremstyle{plain}
\newtheorem{thm}{Theorem}\crefname{thm}{theorem}{theorems}
\crefname{lem}{lemma}{lemmas}
\crefname{prop}{proposition}{propositions}
\crefname{cor}{corollary}{corollaries}

\crefname{obs}{observation}{observations}

\theoremstyle{definition}
\crefname{dfn}{definition}{definitions}

\newtheorem*{assmp*}{Assumption}

\crefname{innerAssmp}{assumption}{assumptions}

\crefname{assmpF}{Assumption}{Assumptions}

\crefname{assmpQ}{Assumption}{Assumptions}

\crefname{assmpA}{Assumption}{Assumptions}

\theoremstyle{remark}

\newtheorem{exmpl}{Example}\crefname{exmp}{example}{examples}

\crefname{figure}{figure}{figures}

\newcounter{exmp}

\title{On forward invariance in Lyapunov stability theorem for local stability
}
\author{Dai {\sc Zusai}\thanks{Department of Economics, Temple University. E-mail: \texttt{ZusaiDPublic@gmail.com}.}}
\date{\today}

\begin{document}
\thispagestyle{empty}
\maketitle
\begin{abstract}
	Forward invariance of a basin of attraction is often overlooked when using a Lyapunov stability theorem to prove local stability; even if the Lyapunov function decreases monotonically in a neighborhood of an equilibrium, the dynamic may escape from this neighborhood. In this note, we fix this gap by finding a smaller neighborhood that is forward invariant. This helps us to prove local stability more naturally without tracking each solution path. Similarly, we prove a transitivity theorem about basins of attractions without requiring forward invariance.

	
\end{abstract}

\section{Introduction}
The idea of Lyapunov stability theorem or Lyapunov's direct method is intuitive: if we find a mapping (\textit{Lyapunov function}) from the current state of a dynamic to a real number such that i) the function attains a local minimum only at an equilibrium (possibly a set) and ii) its value decreases as long as the current state has not reached the equilibrium, then the equilibrium is stable under the dynamic. With this on hand, (we hope that) we do not have to identify a solution path; we just find a Lyapunov function and see how it behaves in the neighborhood–--in particular, the value and first-order derivatives at each point in the state space. So, we typically find a neighborhood where the decrease in the Lyapunov function is guaranteed, which call here a monotone decrease neighborhood, and expect this neighborhood to be a basin of attraction.

However, a basin of attraction must be forward invariant. (This does not matter for global stability, of course.) Precisely, in known versions of Lyapunov stability theorem (e.g. \citealt{Smirnov01DI}), the monotone decrease must be assured to hold on each solution path. Even if we find a monotone decrease neighborhood, a solution path may escape from this neighborhood and eventually the Lyapunov function may not decrease after the escape. This imposes an additional burden of proof, losing an appeal of the theorem to intuition since we eventually need to identify a solution path.
This is overlooked in practice; e.g. \cite{Sandholm10TELocalStbl} and \cite{ZusaiTBRD} on evolutionary dynamics in games, which we fix in this paper.

Similarly, we would expect transitivity of such basins of attractions. That is, if we find a Lyapunov function that decreases in $X_1$ and attains the minimum in $X_2$ and another that decreases in $X_2$ and attains the minimum in $X^\ast$, then we expect $X^\ast$ to be stable in $X_1$. Again, known versions of transitivity theorems as in \cite{Conley78_IsoInvSets} (see also \cite[Theorem 3]{OSTsBRD}) require $X_1$ to be forward invariant and $X_2$ to be forward and also strongly negative (i.e., backward) invariant.\footnote{Strong negative invariance of $X$ means that, if a solution path (starting at time $0$) visits $X$ at any positive time, then it must have started from $X$ at time 0.}

In applications to economics or game theory, we hope to find a Lyapunov function from economic intuition. Under an agent-based dynamic in a game or economic model, an aggregate of agents' possible gains from adjustment of their choices can be used as a candidate for a Lyapunov function once we find a neighborhood where an agent's revision of the choice incurs negative payoff externality to others' gains from further changes, as generally proven by \cite{ZusaiGains}. However, forward invariance needs more mathematical examination of the dynamic system, which may not be appealing to economic intuition. 

In this paper, we reduce the burden of proof by showing that we can construct a forward invariant (smaller) neighborhood from a monotone decrease neighborhood. This fills the gap in applications, as in the papers mentioned above.  
Further, this helps us to establish a transitivity theorem without requiring forward or negative invariance.

We consider a differential inclusion (a set-valued differential equation) and also an equilibrium set, not necessarily a point. This generalization is needed to cover evolutionary dynamics in games, since Nash equilibrium may constitute a (connected) set and also a transition may not be uniquely specified when there are multiple best responses.


\section{Definitions and theorems}
We consider an autonomous differential inclusion $\bcV$ such as 
$$ \dot\x\in \bcV(\x)$$
on a compact metric $A$-dimensional real space $\cX\subset\R^A$ with $A<\infty$. $T\cX$ stands for the tangent space of $\cX$.\footnote{Below the statements of the definitions follow \cite{SandholmPopText}, a canonical reference book on evolutionary dynamics in games.} As a solution concept for the differential inclusion, we adopt a Carath\'{e}odory solution; that is, a solution path $\{\x_t\}_{t\ge 0}$ must be Lipschitz continuity at every $t\ge 0$ and also differentiable with derivative $\dot\x_t\in \bcV(\x_t)$ at almost every $t$.

Let $X^\ast$ be a nonempty closed set. We say $X^\ast$ is {\bf Lyapunov stable} under $\bcV$ if for any open neighborhood $O$ of $X^\ast$ there exists a neighborhood $O'$ of $A$ such that {\it every} solution path $\{\x^t\}_{t\ge 0}$ that starts from $O'$ remains in $O$. $X^\ast$ is {\bf attracting} if there is a neighborhood $O$ of $X^\ast$ such that {\it every} solution that starts in $O$ converges to $X^\ast$; $O$ is called a basin of attraction to $X^\ast$. If it is the entire space $\cX$, then we say $X^\ast$ is globally attracting. $X^\ast$ is {\bf asymptotically stable} if it is Lyapunov stable and attracting; it is globally asymptotically stable if it is Lyapunov stable and globally attracting.

\paragraph{Lyapunov stability theorem.} 

\begin{thm}[Lyapunov stability theorem]\label{thm:Lyapunov_DI}
	Let $X^\ast$ be a non-empty closed set in a compact metric space $\cX$ with tangent space $T\cX$, and $X'$ be a neighborhood of $X^\ast$. Suppose that  continuous function $W:\cX\to\R$ and lower semicontinuous function $\tilde W:\cX\to\R$ satisfy (a) $W(\x)\ge 0$ and $\tilde W(\x)\le 0$ for all $\x\in X'$ and (b) $\cl X'\cap W^{-1}(0)= X'\cap\tilde W^{-1}(0)=X^\ast$. In addition, assume that $W$ is Lipschitz continuous in $\x\in X'$. i) If a differential inclusion $\bcV:\cX\to T\cX$ satisfies\footnote{$D$ denotes differentiation, so $DW(\x)=\tdif{W}{\x}(\x)=[\tpdif{W}{x_1}(\x),\ldots,\tpdif{W}{x_A}(\x)]$.}
	\begin{equation}
	DW(\x)\dot\x \le \tilde W(\x) \qquad \text{for any }\dot\x\in\bcV(\x)  \label{eq:LCondtn_DI}
	\end{equation}
	whenever $W$ is differentiable at $\x\in X'$, then $X^\ast$ is asymptotically stable under $\bcV$. ii) If $X'$ is forward invariant, i.e., every Carath\'{e}odory solution path $\{\x_t\}$ starting from $X'$ at time $0$ remains in $X'$ for all moments of time $t\in[0,\infty)$, then $X'$ is a basin of attraction to $X^\ast$.
\end{thm}
We call $W$ a \textit{Lyapunov function} and $\tilde W$ a \textit{decaying rate function}. Note that we allow for multiplicity of transition vectors, while requiring functions $W$ and $\tilde W$ to be well defined (the uniqueness of the values) as functions of state variable $\x$, independently of the choice of transition vector $\dot\x$ from $\bcV(\x)$.

In a standard Lyapunov stability theorem (e.g. \citet[\S5.5.3]{Robinson_DynSysEd2}) for a differential equation, a decaying rate function $\tilde W$ is not explicitly required while $\dot W$ is assumed to be (strictly) negative until $\x$ reaches the limit set $X^\ast$. The most significant difference is the requirement of lower semicontinuity of $\tilde W$. This assures the existence of a lower bound on the decaying rate $\dot W(\x)\le \bar w<0$ in a hypothetical case in which $\x$ remained out of an arbitrarily small neighborhood of $X^\ast$ for an arbitrarily long period of time. This excludes the possibility that $\x$ would stay there forever and guarantees convergence to  $X^\ast$ (not only Lyapunov stability, i.e., no asymptotic escape from $X^\ast$).

\citet[Theorem 7]{ZusaiTBRD} modifies the Lyapunov stability theorem for a differential inclusion in \citet[Theorem 8.2]{Smirnov01DI}. While the latter is applicable to a singleton of an equilibrium point, the former allows convergence to a set of equilibria.

\Cref{thm:Lyapunov_DI} in this paper relaxes assumptions in \citet[Theorem 7]{ZusaiTBRD}. The previous version imposes a stronger assumption than \eqref{eq:LCondtn_DI}: every Carath\'{e}odory solution $\{\x_t\}$ starting from $X'$ should satisfy 
\begin{equation}
\dot W(\x_t) \le \tilde W(\x_t) \quad \text{for almost all }t\in[0,\infty). \label{eq:LCondtn_DI_tBRD}
\end{equation}
If $X'$ is forward invariant, then \eqref{eq:LCondtn_DI} implies this condition; thus, \cite{ZusaiTBRD} is straightforwardly applied and we can conclude that $X^\ast$ is asymptotically stable and $X'$ is a basin of attraction, as restated in part ii) of \Cref{thm:Lyapunov_DI}.\footnote{Once monotone decrease in $W$ is confirmed for any solution path as in \eqref{eq:LCondtn_DI_tBRD}, convergence to $W^{-1}(0)$ is obtained simply by using Gr\"{o}nwall's inequality.} In part i) of \Cref{thm:Lyapunov_DI} in the current version, we do not require forward invariance of $X'$; a solution trajectory may escape from $X'$ and thus \eqref{eq:LCondtn_DI_tBRD} may not be maintained. Thus, the current version weakens the assumption.

Besides, \eqref{eq:LCondtn_DI} is assumed for every point in the entire space $\cX$ and $\tilde W$ is assumed to be continuous. By checking the places in the proof where the definition of the domain for condition (i) in the theorem\footnote{In the notation of the current version, the condition reads as $W(\x)\ge 0$ and $\tilde W(\x)\le 0$ for all $\x\in\cX$. Thus it corresponds with condition (a) in \Cref{thm:Lyapunov_DI}.} and continuity of $\tilde W$ were used, one can easily find that it is innocuous to replace simplex $\Delta^A\subset\R^A$ with a closed subset $X^\ast$ of a compact metric space $\cX\subset\R^A$ and relax continuity of $\tilde W$  to lower semicontinuity.\footnote{Specifically, $\breve A$ in the proof \citep[p.25]{ZusaiTBRD} should be defined as a subset of $\cl X'$. The continuity of $\tilde W$ was used to assure the existence of the minimum of $\tilde W$ in $\breve A$; for this, lower semicontinuity is sufficient. Then, with the observation that forward invariance of $X'$ implies \eqref{eq:LCondtn_DI_tBRD}, the proof for the previous version applies to part ii) of the current version.} 

\paragraph{Transivitity theorem.} 
\begin{thm}[Transivitity theorem]\label{clm:Lyap_trans}
	Let $X_1\supset X_2 \supset X^\ast$ be three non-empty subsets of a compact metric space $\cX$; assume that $X^\ast$ is closed and $X_1$ is open. Suppose that two Lipschitz continuous functions $W_1,W_2:X_1\to\R$ and two lower semicontinuous functions $\tilde W_1,\tilde W_2:X_1\to\R$ satisfy the following assumptions: for any $\x\in X_1$,
	\begin{enumerate}[a)]
		\item i) $W_1(\x)\ge 0$, ii) $\tilde W_1(\x)\le 0$, and iii) $\cl X_1\cap {W_1}^{-1}(0)=\cl X_1\cap {\tilde W_1}^{-1}(0)=\cl X_2$;
		\item i) $W_2(\x)\ge 0$, ii) $\left[\x\in X_2 ~\Rightarrow~ \tilde W_2(\x)\le 0 \right]$, and iii) $\cl X_2\cap{W_2}^{-1}(0)=\cl X_2\cap{\tilde W_2}^{-1}(0)=X^\ast$;
		\item $\tilde W_1(\x)+\tilde W_2(\x)\le 0$.
	\end{enumerate}
	Furthermore, assume that
	$$ \text{a-iv) } DW_1(\x) \dot\x \le \tilde W_1(\x), \qquad  \text{b-iv) } DW_2(\x) \dot\x \le \tilde W_2(\x) \qquad\text{for any $\dot\x\in\bcV(\x)$},$$
	whenever $W_1$ and $W_2$ are differentiable at $\x\in X_1$. Then, $X^\ast$ is asymptotically stable under $\bcV$. 
\end{thm}

Conditions a) imply that $W_1$ works as a Lyapunov function for local asymptotic stability of $X_2$ in (a subset of) $X_1$ and conditions b) imply that $W_2$ works as a Lyapunov function for local asymptotic stability of $X^\ast$ in (a subset of) $X_2$. So, we may jump to conclude that $X^\ast$ is asymptotically stable in $X_1$. If $X_1$ is indeed forward invariant and thus a basin of attraction to $X_2$ and $X_2$ is forward and also strongly negative invariant, then the transitivity theorem as in \citet[\S II.5.3.D]{Conley78_IsoInvSets} and \citet[Thorem 3]{OSTsBRD} guarantees asymptotic stability of $X^\ast$ with $X_1$ being a basin of attraction to $X^\ast$.

However, we do not assume invariance of these sets in our theorem. Even though \Cref{thm:Lyapunov_DI} eventually assures the existence of some forward invariant subset of $X_2$ from conditions b), the Lyapunov function $W_1$ in conditions a) guarantees convergence only to $X_2$, but not necessarily to this forward invariant subset. Furthermore, the standard transitivity theorem also requires strongly negative invariance of $X_2$, which may not be satisfied by the basin of attraction that we could find using \Cref{thm:Lyapunov_DI}.

The above theorem  avoids this issue by imposing condition c) to hold in the whole $X_1$, which we use to construct a Lyapunov function in $X_1$ to $X^\ast$. Then, we apply \Cref{thm:Lyapunov_DI} and thus the basin of attraction to $X^\ast$ is smaller than $X_1$. When current state $\x$ is in the interim subset $X_2$, conditions a) and b-ii) imply condition c). Condition c) deals with the case that $\x$ is still out of $X_2$ and thus $\tilde W_2(\x)$ may be positive (and thus $W_2$ may be increasing over time). Condition c) requires this to be suppressed by $\tilde W_1$, which must be negative in the entire $X_1$ by condition a-iii).
 
\paragraph{Applications in game theory.} Local stability of an equilibrium is one of the fundamental issues in game theory. A game may exhibit multiple Nash equilibria and thus each equilibrium may not be globally stable; thus, while we investigate global stability of the set of equilibria, we check local stability of each isolated equilibrium or each isolated connected set of equilibria. 

A canonical condition to derive local stability under economically reasonable dynamics is negative payoff externality, specifically called \textit{self-defeating externality }by \citet{HofSand09JET_StableGames}. We regard \textit{a (population) game} as a mapping (payoff function) $\F$ from a distribution of strategies among (continuously many) agents $\x\in \Delta^A$ to a payoff vector $\bpi=\F(\x)\in\R^A$. Self-defeating externality boils down to negative semidefiniteness of $\z\cdot D\F(\x)\z$; a marginal deviation $\z$ in the strategy distribution from $\x$ triggers the change in the payoff vector, which is approximated as $D\F(\x)\z$. Self-defeating externality imposes a negative correlation between $\z$ and $D\F(\x)\z$; a strategy whose share increases by this deviation should face a decrease in its payoff and thus becomes less disadvantageous. 

From this condition on $\F$, economists naturally expect agents to return to the equilibrium, while it needs to formulate how agents revise their choices of strategies in response to payoff changes. An evolutionary dynamic is a dynamic of the strategy distribution, constructed as a mean dynamic of agents each of whom revises its own strategy upon a receipt of a revision opportunity following a Poisson process \citep{SandholmPopText}. \cite{Sandholm10TELocalStbl} considers an equilibrium that is essentially characterized by self-defeating externality, called a regular evolutionary stable state and attempts to prove its local stability under several canonical classes of evolutionary dynamics, such as excess payoff dynamics and pairwise payoff comparison dynamics. Similarly, \cite{ZusaiTBRD} proves it for another class of dynamics, called tempered best response dynamics. These papers refer to a standard version of Lyapunov stability theorem as mentioned after \Cref{thm:Lyapunov_DI}, which requires monotone decrease in the Lyapunov function along with each solution path. However, these papers confirm its decrease only at each point in a neighborhood of a regular ESS, where the self-defeating externality holds. It is left unchecked whether this neighborhood is forward invariant. Our \Cref{thm:Lyapunov_DI} fixes this overlooked point.

In both the two papers, the Lyapunov function is decomposed to two parts; one is the aggregate of possible payoff gains for agents from revisions of strategies and another is the mass of agents who currently choose the strategies that are to be abandoned in the regular ESS. Thanks to this decomposition, we can apply \Cref{clm:Lyap_trans}. \cite{Zusai_dynStblESS} uses it to generalize their results to a broader class of economically natural dynamics.

\section{Proofs}
\subsection*{Proof of \Cref{thm:Lyapunov_DI}}

\begin{proof}
	Here we prove the difference in part i) from \citet[Theorem 7]{ZusaiTBRD}. For this, we focus on the case of $X'\subsetneq \cX$ and find a forward invariant subset of $X'$. Once we find it, any Carath\'{e}odory solution starting from the forward invariant subset remains there and thus satisfies \eqref{eq:LCondtn_DI_tBRD} as \eqref{eq:LCondtn_DI} holds for $\x=\x_t$ at each time $t\in\R_+$. Then, \citet[Theorem 7]{ZusaiTBRD} is applied and assures asymptotic stability of $X^\ast$ while having the forward invariant subset as a basin of attraction.
	
	First, construct a distance from point $\x\in\cX$ to $X^\ast$ based on the metric on $\cX$, say $d:\cX\times\cX\to\R_+$, by
	$$ d_\ast(\x) \coloneqq \min_{\x^\ast\in X^\ast} d(\x,\x^\ast).$$
	Since $X^\ast$ is a non-empty compact set and $d(\x,\x^\ast)$ is continuous in $\x^\ast$ when $\x$ is fixed, Weierstrass theorem assures the existence of the minimum in the above definition of $d_\ast(\x)$. This $d_\ast$ satisfies
	$$ d_\ast(\x)\ge 0; \qquad d_\ast(\x)=0 \ \Leftrightarrow\ \x\in X^\ast.$$
	
	Let $\bar d$ be the shortest distance from the complement of $X'$ to $X^\ast$:
	\begin{equation}
	\bar d \coloneqq \min_{\x\in \cX\setminus X'} d_\ast(\x). \label{eq:prf:Lyapunov_DI_dfnBarD}
	\end{equation}
	Maximum theorem guarantees continuity of $d_\ast:\cX\to\R_+$ by continuity of $d(\x,\x^\ast)$ in both $\x$ and $\x^\ast$. Besides, $\cX\setminus X'$ is a non-empty compact subset by $X'\subsetneq \cX$ and the openness of $X'$. Hence, the minimum in \eqref{eq:prf:Lyapunov_DI_dfnBarD} exists. It follows that 
	\begin{equation}\label{eq:prf:Lyapunov_DI_barD}
	\bar d>0; \qquad d_\ast(\x)<\bar d \ \Rightarrow\ \x\in X'. 
	\end{equation}
	
	\begin{figure}[t!]
		\centering
		\includegraphics[scale=0.6]{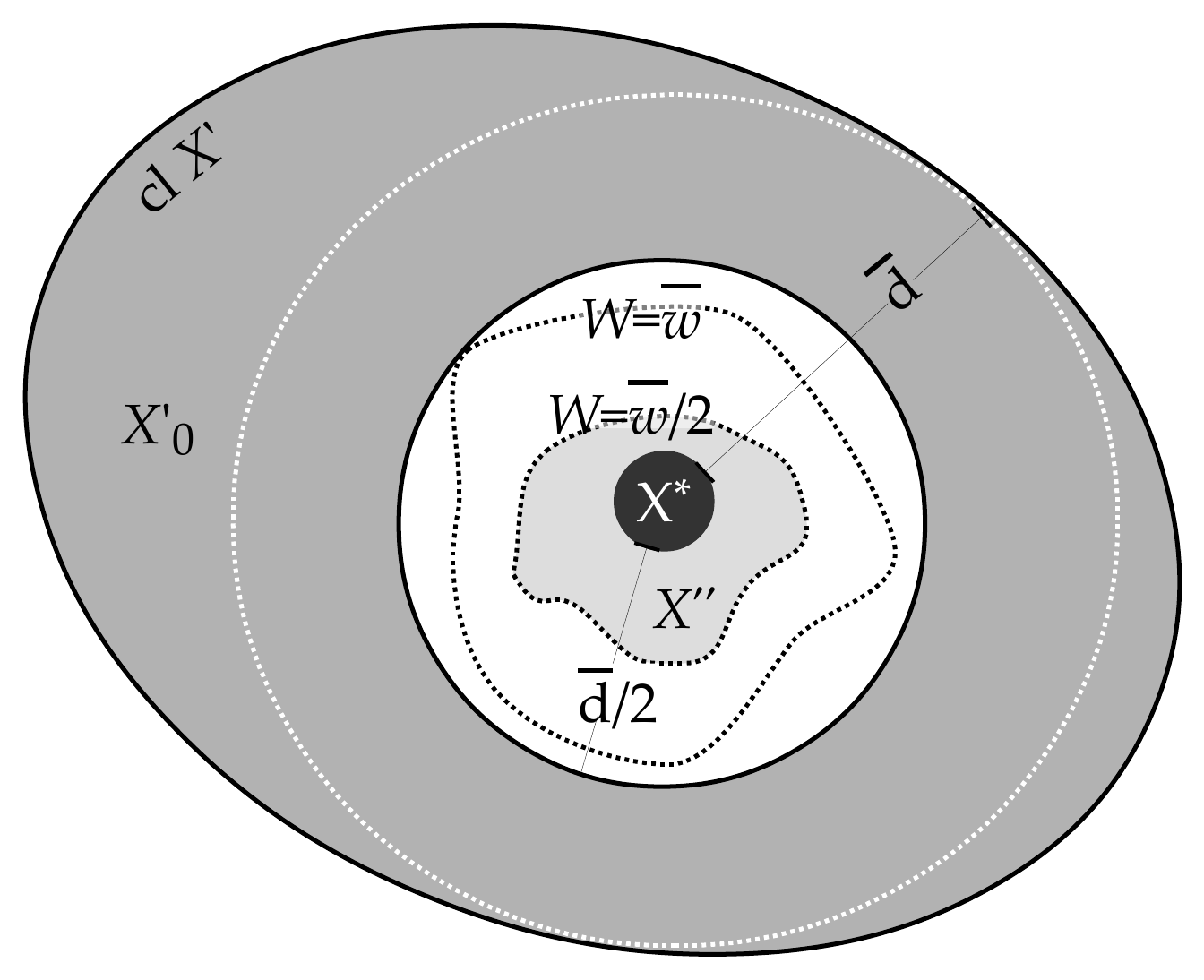}
		\caption{Sets in the proof of \Cref{thm:Lyapunov_DI}. $X^\ast$ is the black area in the center and $X''$ is the light gray area. $\cl X'$ is the entire oval, with the outermost outline. $X'_0$ is the dark gray area, including the both boundaries.}
	\end{figure}
	
	Define set $X'_0\subset\cl X'$ by
	$$ X'_0 \coloneqq \cl X' ~\cap ~ {d_\ast}^{-1}([\bar d/2, \infty)).$$
	Since both $\cl X'$ and ${d_\ast}^{-1}([\bar d/2, \infty))$ are closed, $X'_0$ is closed and thus compact in $\cX$. It is not empty, as proven here. Suppose $X_0'=\emptyset$; then, any $\x\in\cX$ with $d_\ast(\x)\ge \bar d/2$ must be out of $\cl X'$. On the other hand, since $\cl X'$ is not empty, $X'$ has at least one boundary point $\x^0$; then, $d_\ast(\x^0)\ge \bar d$.\footnote{We can make a sequence converging to $\x^0$ from elements out of $X'$, whose distance from $X^\ast$ cannot be smaller than $\bar d$ by \eqref{eq:prf:Lyapunov_DI_dfnBarD}.} By the former statement, this implies $\x^0\notin\cl X'$ but it contradicts with $\x^0$ being on the boundary of $X'$; hence, $X'_0$ cannot be empty.
	
	Let $\bar w$ be the minimum of $W$ in $X'_0$; 
	$$ \bar w \coloneqq \min_{\x\in X'_0} W(\x).$$
	Since $X'_0$ is compact and nonempty and $W$ is (Lipschitz) continuous, the minimum exists. Furthermore, it is positive; we have $X_0'\subset\cl X'$ by construction and $W(\x)\ge 0$ for all $\x\in\cl X'$ by condition (a) and continuity of $W$, while no element $\x\in X'_0$ belongs to $X^\ast$ since $d_\ast(\x)\ge \bar d>0$ for any $\x\in X'_0$. Because $X^\ast =\cl X'\cap W^{-1}(0)$ by condition (b) and $X'_0\subset \cl X'$, it implies $\x\in X'_0 \Rightarrow W(\x)>0$. Hence we have $\bar w>0$; by the definition of $\bar w$, we have
	\begin{equation}\label{eq:prf:Lyapunov_DI_barW}
	\underbrace{\left[ \x\in\cl X' \text{ and } d_\ast(\x) \ge \bar d/2\right]}_{\text{ i.e., }\x\in X'_0}  \ \Rightarrow\ W(\x)\ge \bar w.
	\end{equation}
	
	Define set $X''\subset X'$ by
	\begin{equation}
	X''=W^{-1}([0,\bar w/2)) \cap X'. \label{eq:prf:Lyapunov_DI_dfnX''}
	\end{equation}
	This set is an (open) neighborhood of $X^\ast$ by $X^\ast\subset X''$, since  $W=0$ at anywhere in $X^\ast$ and $X^\ast \subset X'$. Now we prove that $X''$ is wholly contained in set $d_\ast^{-1}([0,\bar d/2))$. Assume that there exists $\x\in X''$ such that $d_\ast(\x)\ge \bar d/2$. These jointly imply $W(\x)\ge \bar w$ by \eqref{eq:prf:Lyapunov_DI_barW} since $\x\in X''\subset X'\subset \cl X'$. However, this contradicts with $W(\x)\in [0,\bar w/2)$ for $\x$ to belong to $X''$. Hence, we have
	\begin{equation}
	\x\in X'' \ \Rightarrow\ d_\ast(\x)< \bar d/2. \label{eq:prf:Lyapunov_DI_X''}
	\end{equation}
	
	Now we prove $X''$ is forward invariant. To verify it by contradiction, assume that there is a Carath\'{e}odory solution trajectory $\{\x_t\}$ starting from $X''$ but escaping $X''$ at some moment of time:
	\begin{equation}\label{eq:prf:Lyapunov_DI_Hypo}
	\x_0\in X'', \qquad \text{ and }\x_T\notin X'' \text{ at some }T>0.
	\end{equation}
	The statement $\x_T\notin X''$ means $\x_T\notin X'$ or $W(\x_T)>\bar w/2$ by \eqref{eq:prf:Lyapunov_DI_dfnX''}. In the former case, we have $d_\ast(\x_T)\ge \bar d$ by \eqref{eq:prf:Lyapunov_DI_dfnBarD} while $d_\ast(\x_0)< \bar d/2$ by \eqref{eq:prf:Lyapunov_DI_X''}. By continuity of $d_\ast(\x)$ in $\x$ and of $\x_t$ in $t$ on a Carath\'{e}odory solution trajectory $\{\x_t\}$, $d_\ast(\x_t)$ is continuous in $t$; hence, there exists a moment of time $T'\in (0,T)$ such that $d_\ast(\x_{T'})=0.9 \bar d \in(0.5 \bar d, \bar d)\subset (d_\ast(\x_0), d_\ast(\x_T))$. At this point, $\x_{T'}\notin X''$ by \eqref{eq:prf:Lyapunov_DI_X''} while $\x_{T'}\in X'$ by \eqref{eq:prf:Lyapunov_DI_barD}; thus, $W(\x_{T'})\ge \bar w/2$ by \eqref{eq:prf:Lyapunov_DI_dfnX''}. Hence, the first case of escaping $X''$ implies the existence of $T'>0$ such that
	$$  W(\x_{T'})\ge \bar w/2 \qquad\textit{ and }\qquad \x_{T'}\in X'. $$
	In the second (but not the first) case, we have $W(\x_T)\ge \bar w/2$ but $\x_T\in X'$; that is, the above statement holds with $T'=T
	.$
	
	This implies the existence of $\bar T\in(0, T']$ such that
	\begin{equation}
	W(\x_{\bar T})\ge \bar w/2, \qquad\text{ and }\qquad \left[ \x_t\in X' \ \text{ for all }t< \bar T \right]. \label{eq:prf:Lyapunov_DI_barT}
	\end{equation}
	To prove it, assume $\x_{t'}\notin X'$ at some $t'< T'$, i.e., the negation of the latter condition with $\bar T=T'$; if there is no such $t'\le T'$, then it suggests that the claim \eqref{eq:prf:Lyapunov_DI_barT} holds at $\bar T=T'$ by the fact  $W(\x_{T'})\ge \bar w/2$.
	%
	%
	By \eqref{eq:prf:Lyapunov_DI_barD}, the hypothesis $\x_{t'}\notin X'$ implies $d_\ast(\x_{t'})\ge \bar d$. Again, by continuity of $d_\ast(\x_t)$ in $t$, the set $\{t\le t' \mid  d_\ast(\x_t)\ge \bar d\}$ is closed and thus compact; by the fact $d_\ast(\x_0)<\bar d/2$, this implies the existence of the minimum $\bar T$ in this set and $\bar T>0$. That is, we have $d_\ast(\x_t)<\bar d$ for all $t< \bar T$ while $d_\ast(\x_{\bar T})=\bar d$. The former implies $\x_t\in X'$ for all $t<\bar T$ by \eqref{eq:prf:Lyapunov_DI_barD} and the latter implies $W(\x_{\bar T})\ge \bar w$ by $\x_{\bar T}=\lim_{t\to\bar T}x_t\in\cl X'$ and \eqref{eq:prf:Lyapunov_DI_barW}. Thus, the above claim \eqref{eq:prf:Lyapunov_DI_barT} holds at this $\bar T\in(0, T']$. Since condition (a) and \eqref{eq:LCondtn_DI} hold almost everywhere in $X'$, we have $\dot W(\x_\tau)\le \tilde W(\x_\tau)\le 0$ at almost all $\tau< \bar T$;\footnote{A Carath\'{e}odory solution trajectory is differentiable at almost all moments of time, though it may not be so at all moments.} thus, we have
	$$ W(\x_{\bar T})\le W(\x_0)+\int_0^{\bar T}  \tilde W(\x_\tau) d\tau\le W(\x_0).$$
	Since $W(\x_0)<\bar w/2$ by $\x_0\in X''$, we have $W(\x_{\bar T})<\bar w/2$ in \eqref{eq:prf:Lyapunov_DI_barT}. This contradicts with $W(\x_{\bar T})\ge \bar w/2.$
	
	Therefore, the hypothesis \eqref{eq:prf:Lyapunov_DI_Hypo} cannot hold: any Carath\'{e}odory solution trajectory $\{\x_t\}$ starting from $X''$ cannot escape $X''$ at any moment of time. That is, $X''$ is forward invariant. 
\end{proof}

\subsection*{Proof of \Cref{clm:Lyap_trans}}
\begin{proof}
	Define a Lyapunov function $W:X_1\to\R$ and a decaying rate function $\tilde W:X_1\to\R$ by 
	$$ W(\x) \coloneqq 2 W_1(\x)+W_2(\x), \qquad  \tilde W(\x) \coloneqq 2 \tilde W_1(\x)+\tilde W_2(\x) \qquad\text{ for each }\x\in X_1.$$
	Lipschitz continuity of $W_1$ and $W_2$ and lower semicontinuity of $\tilde W_1$ and $\tilde W_2$ are succeeded to those of $W$ and $\tilde W$,  respectively.
	It is immediate from assumptions a-i,iv), b-i,iv) and c) to see that 
	\begin{align}
	W(\x)= 2 W_1(\x)+W_2(\x) \ge 0, \notag\\
	\tilde W(\x)= \tilde W_1(\x)+\{\tilde W_1(\x)+\tilde W_2(\x)\}\le 0, \notag\\
	DW(\x)\dot\x =2 DW_1(\x)\dot\x+ DW_2(\x)\dot\x \le 2 \tilde W_1(\x)+\tilde W_2(\x) =\tilde W(\x)
	\label{eq:prf_Lyap_trans_WLyap1}
	\end{align}
	for any $\x\in X_1, \dot\x\in\bcV(\x)$ (for the last equation assuming that $W_1$ and $W_2$ are differentiable at $\x$). 
	
	Further, since $X^\ast\subset X_2$, it follows assumptions a-iii) and b-iii) that $W(\x)=\tilde W(\x)=0$ if $\x\in X^\ast$; thus $X^\ast$ is contained in $\cl X_1\cap W^{-1}(0)$ and $\cl X_1\cap \tilde W^{-1}(0)$ by $X^\ast\subset X_2\subset X_1\subset \cl X_1$. In contrary, assume $W(\x)=0$ at $\x\in\cl X_1$ first. By assumptions a-i) and b-i), it must be the case that $W_1(\x)=0$ and $W_2(\x)=0$. The former implies $\x\in\cl X_2$ by assumption a-iii). Together with this, the latter implies $\x\in X^\ast$ by assumption b-iii). Separately from this, now assume $\tilde W(\x)=0$ at $\x\in\cl X_1$. By assumptions a-ii) and c), it must be the case that $\tilde W_1(\x)=0$ and $\tilde W_1(\x)+\tilde W_2(\x)=0$.\footnote{Note that the latter condition alone cannot assure $\tilde W_2(\x)=0$, since $\tilde W_2(\x)$ could take a positive value unless $\x$ is in $X_2$.} The former implies $\x\in\cl X_2$ by assumption a-iii); besides, by plugging the former into the latter, we have $\tilde W_2(\x)=0$. These two statements jointly imply $\x\in X^\ast$ by assumption b-iii). In sum, we have verified 
	\begin{equation}
	\cl X_1\cap W^{-1}(0)=\cl X_1\cap {\tilde W}^{-1}(0)=X^\ast. \label{eq:prf_Lyap_trans_WLyap2}
	\end{equation}
	Note that the first equality is due to the fact that $X^\ast\subset X_1$ and thus $X^\ast \cap \bd X_1=\emptyset$ since $X_1$ is open. 
	
	We have verified all the assumptions in \Cref{thm:Lyapunov_DI}; therefore, $X^\ast$ is asymptotically stable. Notice that $X_1$ may not be forward invariant, but part i) of \Cref{thm:Lyapunov_DI} assures that we can make some subset of $X_1$ as a basin of attraction to $X^\ast$.
\end{proof}

\renewcommand\refname{References}
\bibliography{../DZbib}

\end{document}